\documentclass[reqno, 12pt]{amsart}
\usepackage{amssymb,amsmath,amsfonts,amsthm,comment,mathrsfs,times,graphicx}
\usepackage[bookmarksnumbered, colorlinks, plainpages]{hyperref}
\usepackage{color}
\usepackage[english]{babel}
\usepackage[all,cmtip]{xy}
\usepackage{lmodern}
\usepackage{enumitem}
\usepackage{geometry}
\newcounter{dummy} \numberwithin{dummy}{section}
\newtheorem{theo}[dummy]{Theorem }
 
 \newtheorem{lem}[dummy]{Lemma}
 \newtheorem{pr}[dummy]{Proposition}
 
\newtheorem{ef}[dummy]{Definition}

\newtheorem{rem}[dummy]{Remark}

\newtheorem{conj}{Conjecture}
\title[ Half-integral weight modular forms and real quadratic $p$-rational fields
]{Half-integral weight modular forms and real quadratic $p$-rational
fields}
\author[Jilali Assim, Zakariae Bouazzaoui ]{Jilali Assim$^{(1)}$, Zakariae Bouazzaoui$^{(2)}$}
\address{$^{(1)}$ Moulay Ismail University,
 Department of mathematics,
Faculty of sciences, Meknès, B.P. 11201 Zitoune, Meknès, Morocco.}
\email{\textcolor[rgb]{0.00,0.00,1.00}{j.assim@yahoo.fr}}
\address{$^{(2)}$ Moulay Ismail University,
 Department of mathematics,
Faculty of sciences, Meknès, B.P. 11201 Zitoune, Meknès, Morocco.}
\email{\textcolor[rgb]{0.00,0.00,1.00}{z.bouazzaoui@edu.umi.ac.ma}}
\begin{document}
\maketitle
\renewcommand{\abstractname}{Abstract}
\begin{abstract}
Using half-integral weight modular forms we give a criterion for the
existence of real quadratic $p$-rational fields. For $p=5$ we prove
the existence of infinitely many real quadratic $p$-rational fields.
\end{abstract}

\section{Introduction}

The Dedekind zeta function of an algebraic number field encodes a
lot of arithmetic information of the field. For a number field $F$,
let $\mathcal{O}_F$ denote the ring of its integers. For each
integer $m$, let $\zeta_{F}^{*}(m)$ denote the leading non-zero
coefficient in the Taylor expansion of the Dedekind zeta function of
$F$. Dirichlet's class number formula reads:
\begin{equation}
\zeta_{F}^{*}(0)=-\frac{h_F}{w_F}.R_F,
\end{equation}
where $h_F$ is the class number of $F$, $w_F$ is the number of roots
of unity in $F$ and $R_F$ is the Dirichlet regulator. We are
interested with the divisibility, by odd prime numbers $p$, of the
special values of Dedekind zeta functions of real quadratic fields
at odd negative integers, these values are closely related to the
orders of certain cohomology groups.\\
Let $S$ be a finite set of primes. Denote by $F_S$ the maximal
pro-$p$-extension of $F$ which is unramified outside $S$ and let
$G_{S}(F)$ be its Galois group. The field $F$ is called $p$-rational
if the Galois group $G_{S_p}(F)$ of the extension $F_{S_p}/F$ is
pro-$p$-free (with rank $1+r_2$, $r_2$ being the number of complex
primes of $F$), where $S_p$ is the set of primes of $F$ above $p$.
If $F$ is totally real, we prove in section $2$ that $F$ is
$p$-rational precisely when $v_p(\zeta_F(2-p))=-1$, where $v_p$
denotes the $p$-adic valuation. We use this characterization to
study $p$-rationality of real quadratic fields. The notion of
$p$-rational fields has been introduced to construct extensions of
$\mathbf{Q}$ satisfying  the Leopoldt conjecture
\cite{Movahhedi-Nguyen}. Recently, R. Greenberg \cite{Greenberg}
used $p$-rational number fields to construct (in a non geometric
manner) Galois representations with open image in
$\mathrm{GL}_{n}(\mathbf{Z}_p)$ for $n\geq3$. This paper is
motivated by the study of $p$-rationality of multi-quadratic number
fields, for which Greenberg  formulated the following conjecture:

\begin{conj}\label{p-rational GC}  (\cite[Conjecture 4.2.1]{Greenberg})
For any odd prime number $p$ and any integer $t\geq1$, there is a
$p$-rational field $F$ such that
$\mathrm{Gal}(F/\mathbf{Q})\cong(\mathbf{Z}/2\mathbf{Z})^t$.
\end{conj}

The conjecture is true for $t=1$, since for every odd prime number
$p$, there is infinitely many $p$-rational imaginary quadratic
fields (cf. \cite[Proposition 4.1.1]{Greenberg}). The case $t\geq2$
leads to the study of $p$-rationality of real quadratic fields,
which is the aim of this paper. After relating the $p$-rationality
to special values of $L$-functions, we use the theory of modular
forms to obtain our results. Roughly speaking, we use
Cohen-Eisenstein series \cite{Cohen}, which are modular forms of
half integer weight, and whose Fourier coefficients involve special
values of $L$-functions of quadratic fields. Multiplying such
modular forms by theta series produces integer weight modular forms,
and the resulting Fourier coefficients are studied to deduce
divisibility properties of values of $L$-functions. As a consequence
we give for $p=5$ the existence of infinitely many real quadratic
$5$-rational fields, a similar result for $p=3$ was given implicitly
by D. Byeon in \cite{Byeon} using the same techniques.

\begin{theo}\label{the prime 5}
There are infinitely many fundamental discriminants $d>0$ such that
$\mathbf{Q}(\sqrt{d})$ is $5$-rational.
\end{theo}

The study of $p$-rationality of real quadratic fields is more subtle
than the study of $p$-rationality of imaginary quadratic fields,
because of complications due to the existence of non-trivial units.
Using Cohen-Eisenstein series, Theorem \ref{wieferich} below gives a
sufficient condition for the existence of a real quadratic
$p$-rational field, with some arithmetic properties, for every prime
number $p\geq5$. More precisely, let $f=\sum_{n\geq0}a(n)q^n$ be an
integer weight modular form for the congruence subgroup $\Gamma(N)$,
$N\geq1$, with coefficients in the ring of integers of a number
field. By a result of Serre \cite[page 20-19]{Serre76}, there is a
set of primes $\ell\equiv1\pmod{Np^2}$ of positive density for which
\begin{equation}\label{Serre congruence}
f|T(\ell)\equiv2f\pmod{p^2},
\end{equation}
where $T(\ell)$ denotes the Hecke operator associated to the prime
number $\ell$ \cite[page.153]{Koblitz}.

Let $\mathcal{L}=\{\ell_1,...,\ell_s\}$ be a finite set of odd
primes. For every positive square free integer $t$, let $f$ be an
element of the space $M_p (\Gamma_1
(4p^{2}t\prod_{i=1}^{s}\ell_{i}^{4}))$, obtained by multiplication
of half integer weight modular forms (Cohen-Eisenstein series and
theta series). Denote by $\mathcal{S}_t$ the set of primes $\ell$
satisfying (\ref{Serre congruence}) for $f$. We
make the following hypothesis:\\

$(H_p)$:\; There exist a square free integer $t$ and a prime number
$\ell\in\mathcal{S}_t$ such that $\ell=ta^{2}+b^{2}$
\quad\quad and $b$ is a prime number for which $p$ is non-Wieferich.\\

\begin{theo}\label{wieferich}
Let $\mathcal{L}=\{\ell_1,...,\ell_s\}$ be a finite set of odd
primes. Let $p\geq5$ be a prime number. Assume that hypotheses
$(H_p)$ is satisfied for some prime number $\ell$. Then there is a
real quadratic $p$-rational field $\mathbb{Q}(\sqrt{d})$ for some
fundamental discriminant $d<\ell$ such that $(\frac{d}{\ell_k})=1$
for every $\ell_k\in\mathcal{L}$, where $(\frac{.}{\ell_k})$ denotes
the Legendre symbol.
\end{theo}

\section{$p$-rationality of quadratic fields}

Let $p$ be an odd prime number and let
$\mathcal{O}'_F=\mathcal{O}_F[\frac{1}{p}]$ be the ring of
$p$-integers of $F$, then the field $F$ is called $p$-rational if
the étale cohomology group
$H^{2}(\mathcal{O}'_F,\mathbb{Z}/p\mathbb{Z})$ vanishies
\cite{Movahhedi-Nguyen}, \cite{Movahhedi}, \cite{Movahhedi90}. In
general, for every integer $i$, if
$H^{2}(\mathcal{O}'_F,\mathbb{Z}/p\mathbb{Z}(i))=0$ then we say that
the field $F$ is $(p,i)$-regular \cite{Assim}. If $F$ is totally
real, the information about the $p$-rationality and the
$(p,i)$-regularity of $F$ are contained in special values of the
Dedekind zeta function $\zeta_F$ at odd negative integers. More
precisely, as a consequence of the Main Conjecture in Iwasawa theory
for totally real number fields and odd primes $p$ proved by A.Wiles,
we obtain the following case of Lichtenbaum conjecture: for any even
positive integer $i\geq2$, and any totally real number field $F$, we
have
\begin{equation}\label{formula (2)}
w_i(F)\zeta_F(1-i)\sim_p |H^2(\mathcal{O}'_F,\mathbb{Z}_p(i))|,
\end{equation}
where $w_i(F)$ is the order of the group
$H^0(F,\mathbb{Q}_p/\mathbb{Z}_p(i))$ and $\sim_p$ means that they
have the same $p$-adic valuation. Moreover, a periodicity statement
on cohomology groups gives that
$$H^2(\mathcal{O}'_F,\mathbb{Z}/p\mathbb{Z}(i))\cong H^2(\mathcal{O}'_F,\mathbb{Z}/p\mathbb{Z}(j)),$$
whenever $i\equiv j\pmod{p-1}$. Then we have the following
proposition:

\begin{pr}
Assume that $i\geq2$ is an even integer, then a totally real number
field $F$ is $(p,i)$-regular if and only if $w_{i}(F)\zeta_F(i)$ is
a $p$-adic unit.
\end{pr}

\noindent{\textbf{Proof}.} The proof follows from (\ref{formula
(2)}) and Proposition $2.4$ of \cite{Assim}.$\hfill\blacksquare$
\vskip 6pt

Suppose that $F$ is a totally real number field of degree $g$. Let
$v_p$ be the $p$-adic valuation. We have for even positive integers
$i$ the following result \cite[Theorem 6]{Serre71}:

\begin{theo}
Let $p$ be an odd prime number.
\begin{enumerate}
    \item if $gi\equiv0\pmod{p-1}$,
    $v_p(\zeta_{F}(1-i))\geq-1-v_p(g)$;
    \item if $gi\not\equiv0\pmod{p-1}$, $v_p(\zeta_F(1-i))\geq0$.
\end{enumerate}
\end{theo}

In particular, we have
\begin{equation}
v_p(\zeta_{F}(2-p))\geq-1-v_p(g).
\end{equation}
In \cite[Section 3.7]{Serre}, it is suggested that often
$v_p(\zeta_{F}(2-p))\leq-1$. Using Formula (\ref{formula (2)}) we
relate the $p$-rationality of $F$ to the
special value $\zeta_F(2-p)$ in the following way:\\
Let $p\geq3$ be a prime which is unramified in $F$, then
\begin{equation}\label{Main equivalence}
F\; \hbox{is}\;
p\hbox{-rational}\;\Leftrightarrow\;\;v_p(\zeta_F(2-p))=-1.
\end{equation}
For $F=\mathbf{Q}(\sqrt{d})$ a real quadratic fields, the Dedekind
zeta function of $F$ satisfies
\begin{equation*}
\zeta_{F}(2-p)=\zeta_{\mathbf{Q}}(2-p)L(2-p,\chi_d).
\end{equation*}
Since the field $\mathbf{Q}$ is $p$-rational for every odd prime
number $p$ (which is equivalent to say that
$v_p(\zeta_{\mathbf{Q}}(2-p))=-1$), we have the following
proposition:

\begin{pr}\label{propo 2.1}
Assume that $p\nmid d$, then the field $\mathbf{Q}(\sqrt{d})$ is
$p$-rational precisely when $v_p(L(2-p,\chi_d))=0$.
\end{pr}

This is the motivation behind using the half-integer weight modular
forms called Cohen-Eisenstein series described in the next section.

\section{Cohen-Eisenstein series}

Let $d<0$ be a fundamental discriminant and denote by $h(d)$ the
class number of $\mathbf{Q}(\sqrt{d})$. For a rational prime
$p\geq5$, it is known that if $p\nmid h(d)$ then the field
$\mathbf{Q}(\sqrt{d})$ is $p$-rational. An object which generate
class numbers of imaginary quadratic fields is the $3$-power of the
standard theta series $\theta$ given by the $q$-expansion
$$\theta(q)=1+2q+2q^{4}+2q^9+\cdot\cdot\cdot\cdot$$ This series has been used to prove the
existence of infinitely many $p$-rational imaginary quadratic
fields. More precisely, the series $\theta^3$ is a modular form of
weight $\frac{3}{2}$ for the congruence subgroup $\Gamma_0(4)$.
Write
$$\theta^3(q)=\sum_{n\geq0}r_3(n)q^n,$$ then the coefficient $r_3(n)$ is
the number of times we can write $n$ as a sum of three squares.
These coefficients satisfy, by a Theorem of Gauss, the following
property:\\ If $n>3$ is square free, then
$$r_3(n)=\left\{
           \begin{array}{ll}
             12h(-4n), & \hbox{ $n\equiv1,2,5,6\pmod{8}$;}\\
             12h(-n), & \hbox{$n\equiv3\pmod{8}$.}
           \end{array}
         \right.
$$
By studying divisibility properties by a prime $p$ of the Fourier
coefficients of $\theta^3$, one can deduce the existence of
infinitely many imaginary quadratic fields $\mathbf{Q}(\sqrt{d})$
with $p\nmid h(d)$ \cite{Horie}, \cite{Bruinier}, hence we have
infinitely many $p$-rational imaginary quadratic fields.\\
To study $p$-rationality of real quadratic fields, we shall use the
same approach as for the imaginary quadratic case by considering the
following half-integer weight modular forms.\\

H. Cohen \cite{Cohen} constructed a modular form $\mathcal{H}_{i}$
of weight $i+\frac{1}{2}$ for the congruence subgroup $\Gamma_0(4)$
given by the $q$-expansion:
\begin{equation*}
\mathcal{H}_{i}(q)=\sum_{n\geq0}h(i,n)q^n,
\end{equation*}
such that the Fourier coefficients are essentially given in terms of
special values of Dirichlet $L$-functions of quadratic fields:
$$h(i,n)=\left\{
    \begin{array}{ll}
      \zeta(1-2i), & \hbox{$n=0$;} \\
      0, & \hbox{$(-1)^{i}n\equiv2,3\pmod{4}$;} \\
      L(1-i,\chi_{(-1)^{i}d})\sum_{r|m}\mu(r)\chi_{(-1)^{i}d}(r)r^{i-1}\sigma_{2i-1}(m/r), & \hbox{$(-1)^{i}n=dm^2$.}
    \end{array}
  \right.
$$
where $\mu$ is the M\"{o}bius function, and $\sigma_{s}(n)$ is the
sum of $s$-th powers of the positive divisors of $n$, i.e,
$\sigma_{s}(n)=\sum_{r|n}r^{s}$.\\
Let
$$\epsilon_{d}*\sigma_{2i-1}(m)=\sum_{r|m}\mu(r)\chi_{(-1)^{i}d}(r)r^{i-1}\sigma_{2i-1}(m/r)$$ as a convolution
product of multiplicative functions, where
$$\epsilon_d(r)=\mu(r)\chi_{(-1)^{i}d}(r)r^{i-1}.$$
If $d=1$ we write $\epsilon_d(n)=\epsilon(n)$; this is a
multiplicative function so that
$$h(i,n^{2})=\zeta(1-i)\epsilon*\sigma_{2i-1}(n).$$
Using Proposition \ref{propo 2.1} we see that information about the
$p$-rationality of real quadratic fields are encoded in the Fourier
coefficients of $\mathcal{H}_{p-1}$. In fact, a field
$\mathbf{Q}(\sqrt{d})$ is $p$-rational precisely when
$$v_p(h(p-1,d))=0.$$
Based on this observation, the aim is to evaluate the $p$-adic
valuation of the coefficients of $\mathcal{H}_{p-1}$. For $p=5$, we
shall use the following result of Cohen:

\begin{pr}(\cite[Proposition 5.1]{Cohen})\label{propo4}
Let $r$ be a positive integer. Let $D\equiv0$ or $1\pmod{4}$ be an
integer such that $(-1)^{r}D=|D|$. Then for $r\geq2$
\begin{equation*}
\sum_{n\geq0}(\sum_{s}h(r,\frac{4n-s^2}{|D|}))q^n\in
M_{r+1}(\Gamma_0(D),\chi_D),
\end{equation*} where $M_{r+1}(\Gamma_0(D),\chi_D)$ is the space of
modular forms of weight $r+1$ and level $D$ with character $\chi_D$.
\end{pr}

As a consequence, Cohen \cite[examples, p.277]{Cohen} obtained
formulas such as
$$\sum_{s}h(2,N-s^{2})=\frac{-1}{30}\sum_{r\mid N}(r^{2}+(N/r)^{2})(\frac{-4}{r}).$$
This equality is used by Beyon \cite{Byeon} to prove the existence
of $3$-rational real quadratic fields. Using this approach we prove
Proposition \ref{theorem2} below.\\

For $p\geq5$, we use operators introduced by Shimura \cite{Shimura}
to produce from $\mathcal{H}_{p-1}$ a modular form $G$ with specific
Fourier coefficients. For this, let $f=\sum_{n\geq0}a(n)q^n$ be a
modular form of weight $i+\frac{1}{2}$ in the space
$M_{i+\frac{1}{2}}(N,\chi)$. Let $m>0$ be an integer. The operator
$B_m$ applies $f$ to the modular form
$$f|B_m=\sum_{n\geq0}a(nm)q^{nm},$$ which is an element of
the space $M_{i+\frac{1}{2}}(Nm^2,\chi)$ \cite[Section 3]{Bruinier}.\\
The twist operator is defined as follows. Let $\psi$ be a primitive
Dirichlet character modulo $m$, then twisting $f$ by $\psi$ gives
the modular form
$$f_{\psi}=\sum_{n\geq1}\psi(n)a(n)q^n,$$ which belongs to the space
$M_{i+\frac{1}{2}}(Nm^2,\chi\psi^2)$ \cite[Lemma 3.6]{Shimura}.\\
Combining these two operators we obtain the following modular form
\begin{equation}\label{operators of Shi}
\frac{1}{2}(f-f|B_m)+\frac{1}{2}(f-f|B_m)_{\psi}=\sum_{\psi(n)=1}a(n)q^n.
\end{equation}
Let $\mathcal{H}$ be the modular form defined by
\begin{equation}
\mathcal{H}(q)=p(\mathcal{H}_{p-1}-\mathcal{H}_{p-1}\mid
B_p)(q)=\sum_{n\geq1}h_1(p-1,n)q^n,
\end{equation}
where $h_1(p-1,n)=ph(p-1,n)$. Let
$\mathcal{L}=\{\ell_1,...,\ell_s\}$ be a set of odd prime numbers
and put
$$\mathcal{N}=\{n:\;(\frac{n}{\ell_k})=1\;\;\forall\ell_{k}\in\mathcal{L}\}.$$
Using (\ref{operators of Shi}) for the modular form $\mathcal{H}(q)$
with $m=\ell_1$ and $\psi=(\frac{.}{\ell_1})$, we obtain a modular
form $$G_1(q)=\sum_{(\frac{n}{\ell_1})=1}h_1(p-1,n)q^n.$$ Now do the
same think for $G_1$ and $\ell_2$ to obtain a modular form $G_2$. By
induction on the set $\mathcal{L}$ we obtain the following lemma:

\begin{lem}\label{cohen form}
The modular form
\begin{equation*}
G(q)=\sum_{n\in\mathcal{N}}h_1(p-1,n)q^n
\end{equation*}
is an element of the space
$M_{p-1+\frac{1}{2}}(\Gamma_{1}(4p^{2}\prod_{j=1}^{s}\ell_{j}^{4}))$.
\end{lem}

This modular (with level $\frac{p}{2}$) has been used by Ken Ono in
\cite{Ono} to study class number of real quadratic fields.

\begin{rem}
\begin{enumerate}
\item Let $n$ be a positive integer such that $\gcd{(n,\ell_k)}=1$, for
every $\ell_k\in\mathcal{L}$. Then $n^{2}\in\mathcal{N}$ and the
coefficient $h_1(p-1,n^{2})=\zeta(2-p)\epsilon*\sigma_{2p-3}(n)$ is
a $p$-adic unit, since $\zeta(2-p)$ has $p$-adic valuation $-1$ and
for every prime number $\ell$,
$\epsilon*\sigma_{2i-1}(\ell^n)\equiv1\pmod{p}$. Hence
$\mathcal{H}(q)$ is not a constant modulo $p$.
\item For every integer $n$ such that $h_1(p-1,n)\neq0$, we
have $p\nmid n$.
\end{enumerate}
\end{rem}

\section{ Proof of Theorem \ref{the prime 5}}

In this section, we evaluate sums of Fourier coefficients of the
Cohen-Esenstein series $\mathcal{H}_{p-1}$ for $p=5$, to prove the
existence of infinitely many $5$-rational real quadratic fields.

\begin{pr}\label{theorem2}
Let $\ell,\ell'$ be distinct odd prime numbers such that
$\ell\equiv1\pmod{4}$, $\ell'\equiv3\pmod{4}$ and
$v_5(1-\ell'^4)=1$. Then there is a positive fundamental
discriminant $d_{\ell\ell'}<2\ell\ell'$ such that
$2\ell\ell'=x^{2}+d_{\ell\ell'}y^2$ for some integers $x,y$ and
$\mathbb{Q}(\sqrt{d_{\ell\ell'}})$ is a $5$-rational number field.
\end{pr}

\noindent{\textbf{Proof}.} As a consequence of Proposition
\ref{propo4}, Cohen \cite[page 277]{Cohen} obtained the equality
\begin{equation*}
\sum_{\begin{array}{c}
        s\in\mathbf{Z} \\
        n-s^2\geq0
      \end{array}
}h(5-1,n-s^2)=\frac{1}{300}\sum_{r\mid
n}(r^4+(2\ell/r)^4)(\frac{-4}{r})+\frac{1}{400}\sum_{n=x^2+y^2}(x^4-6x^2y^2+y^4).
\end{equation*}
If $n=2\ell\ell'$, where $\ell'\equiv3\pmod{4}$, then the equation
$2\ell\ell'=x^2+y^2$ has no solutions in integers $x$ and $y$ (it is
well known that an integer $n$ is a sum of two squares if and only
if all primes $\equiv3\pmod{4}$ dividing $n$ have even exponents).
Then, in one hand, we have
$$\frac{1}{400}\sum_{2\ell\ell'=x^2+y^2}(x^4-6x^2y^2+y^4)=0.$$
In the other hand, the sum
$$\sum_{s}h(5-1,n-s^2)$$ is non-empty, since for every integer $m$, the $m$-th coefficient of a Cohen-Eisenstein series
$\mathcal{H}_i$ is non-trivial if and only if
$(-1)^{i}m\equiv0,1\pmod{4}$, and for odd integer $x$ such that
$2\ell\ell'-x^2>0$, we have the congruence
$2\ell\ell'-x^2\equiv1\pmod{4}$.\\
 Moreover, since
$(\frac{-4}{2})=0$, we have the equality:

\begin{equation}
\begin{array}{ccc}
  \sum_{r\mid
2\ell\ell'}(r^4+(2\ell\ell'/r)^4)(\frac{-4}{r}) & = &(1+(2\ell\ell')^4)+(\ell^{4}+(2\ell')^4)-(\ell'^{4}+(2\ell)^4)-((\ell\ell')^{4}+2^4)\\
                                           & = &1-2^{4}-(\ell\ell')^{4}(1-2^4)+\ell^{4}(1-2^4)-\ell'^{4}(1-2^4)\\
                                           &=&(1-2^4)(1-(\ell\ell')^{4}+\ell^{4}-\ell'^4)\\
                                           & = &(1+\ell^4)(1-2^4)(1-\ell'^4).
\end{array}
\end{equation}

Remark that $v_5(1+\ell^4)=0$, $v_5(1-2^4)=v_5(1-\ell'^4)=1$ and
$v_5(300)=2$, hence
$$\frac{1}{300}(1+\ell^4)(1-2^4)(1-\ell'^4)\not\equiv0\pmod{5},$$ which
gives that
$$\sum_{s}h(5-1,6\ell-x^2)\not\equiv0\pmod{5}.$$ Then there
is at least one odd integer $x$ such that
$$v_{5}(h(5-1,2\ell\ell'-x^2))=0.$$ Since $2\ell\ell'$ is not a sum
of two squares, the integer $2\ell\ell'-x^2$ is of the form
$d_{\ell\ell'}y^2$, where $d_{\ell\ell'}>0$ is a fundamental
discriminant. Recall that for every positive integer $i$ and every
fundamental discriminant $(-1)^{i}d$ we have
$$h(i,dn^2)=L(1-i,\chi_{(-1)^{i}d}).\epsilon_{d}*\sigma_{2k-1}(n),$$ and that for a fundamental discriminant $d>0$
and an odd primes $p$ the value  $L(2-p,\chi_{d})$ is of
non-negative $p$-adic valuation \cite[théorème 6]{Serre71}. Then the
fundamental discriminant $d_{\ell\ell'}$ satisfies:
$$v_5(L(2-5,\chi_{d_{\ell\ell'}}))=0.$$ Using the equivalence (\ref{Main equivalence}) this means that the field
$\mathbf{Q}(\sqrt{d_{\ell\ell'}})$ is
$5$-rational.$\quad\hfill\blacksquare$\vskip 6pt

For the proof of Theorem \ref{the prime 5}, we use the following
Lemma:

\begin{lem}
Let $d_1,..,d_m$ be positive fundamental discriminants. Then there
is a positive density of prime numbers $\ell$ for which the
following properties are satisfied:
\begin{enumerate}
    \item $\ell\equiv3\pmod{4}$,
    \item $v_5(1-\ell)=1$,
    \item $(\frac{d_i}{\ell})=1$ for every $i\in\{1,...,m\}$.
\end{enumerate}
\end{lem}

\noindent{\textbf{Proof}.} Let consider the following number field
$F=\mathbf{Q}(\sqrt{d_1},...,\sqrt{d_m},\sqrt{-1},\mu_5,\mu_{25})$,
where as usual, $\mu_n$ denotes the group of $n$-th root of unity.
Since the Galois group
$H=\mathrm{Gal}(F/\mathbf{Q}(\sqrt{d_1},...,\sqrt{d_m}))$ is cyclic,
the \v{C}ebotarev's density theorem gives the existence of prime
$\ell$ (in fact an infinite numbers of primes) such that the
Frobenuis $\mathrm{Frob}_\ell$ at the prime $\ell$ generates the
group $H$. In particular these primes $\ell$ satisfy the properties
$(1)$, $(2)$ and $(3)$ of the lemma.$\quad\hfill\blacksquare$\vskip
6pt

Now we are able to prove Theorem \ref{the prime 5}, which gives the
existence of infinitely many $5$-rational real quadratic fields.

\noindent{\textbf{Proof of Theorem} \ref{the prime 5}.} Suppose that
the set $\mathcal{A}_{5-\mathrm{rational}}$ of positive fundamental
discriminants $d$ such that $\mathbf{Q}(\sqrt{d})$ is $5$-rational
is finite. Put $\mathcal{A}_{5-\mathrm{rational}}=\{d_1,...,d_m\}$.
Let $\ell$ be a prime number satisfying the assumptions
\begin{enumerate}
   \item $\ell\equiv3\pmod{4}$,
   \item $v_5(1-\ell)=1$,
   \item $(\frac{d_i}{\ell})=1$ for every $i\in\{1,...,m\}$.
\end{enumerate}
By Proposition \ref{theorem2}, for every prime number
$\ell'\equiv1\pmod{4}$, there exists a fundamental discriminant
$d_{\ell\ell'}<2\ell\ell'$ such that
$\mathbf{Q}(\sqrt{d_{\ell\ell'}})$ is $5$-rational, hence there
exists $d_j\in\mathcal{A}_{5-rational}$ such that
$d_j=d_{\ell\ell'}$. Moreover, the fundamental discriminant
$d_{\ell\ell'}$ satisfies the equation
$2\ell\ell'=x^2+d_{\ell\ell'}y^2$ for some integers $x,y$, which
implies that $(\frac{d_{\ell\ell'}}{\ell})=(\frac{d_{j}}{\ell})=-1$.
This is a contradiction since the rational prime $\ell$ satisfy
$(\frac{d_k}{\ell})=1$ for every $k\in\{1,...,m\}$.
$\quad\blacksquare$\vskip 6pt

\section{ Proof of Theorem \ref{wieferich}}

In this section, we give a sufficient condition (hypothesis ($H'_p$)
below) for the existence of real quadratic $p$-rational fields for
every prime $p\geq5$. To do this, we shall use the modular form
$$G(q)=\sum_{n\in\mathcal{N}}h_1(p-1,n)q^n$$ in Lemma \ref{cohen
form}, together with the theta series
\begin{equation*}
\theta_t=1+2\sum_{n\geq1}q^{tn^2},
\end{equation*}
where $t>1$ is a fundamental disrciminant (cf. \cite[Section
2]{Shimura} and \cite[Lemma 2]{Serre-Stark}). Recall that $G$ is an
element of the space
$M_{p-1+\frac{1}{2}}(4p^{2}\prod_{j=1}^{s}\ell_{j}^4,\chi_0)$, which
is a non-constant modular form, where $\mathcal{N}=\{n :
(\frac{n}{\ell_k})=1,\;\forall\ell_k \in\mathcal{L}\}$,
$\mathcal{L}=\{\ell_1,...,\ell_s\}$ is a finite set of odd primes
and $\chi_0$ is the trivial character. The series $\theta_t$ belong
to the space $M_{\frac{1}{2}}(4t,\chi_t)$, where $\chi_t$ is the
quadratic character associated to the field $\mathbf{Q}(\sqrt{t})$.
Then the product $G\theta_t$ gives a modular form of weight $p$ in
$M_{p}(4tp^{2}\prod_{j=1}^{s}\ell_{j}^4,\psi)$, where $\psi$ is a
Dirichlet character modulo $4tp^{2}\prod_{j=1}^{s}\ell_{j}^4$. We
write
\begin{equation*}
G\theta_t=\sum_{n\geq0}c(n)q^n,
\end{equation*}
where
\begin{equation*}
c(n)=\sum_{\begin{array}{c}
             n=tx^2+y \\
             y\in\mathcal{N}
           \end{array}
}\alpha(x)h_1(p-1,y).
\end{equation*}
Here $\alpha(x)$ is the $tx^2$-th coefficient of $\theta_t$.

The result of Serre \cite[page 20]{Serre76}, mentioned in the
introduction, gives the existence of a set $\mathcal{S}_t$ of primes
$\ell\equiv1\pmod{4tp^2\prod_{j=1}^{s}\ell_{j}^4}$ of positive
density such that the congruence (\ref{Serre congruence}) holds for
$G\theta_t$, i.e., for every prime number $\ell\in\mathcal{S}_t$ we
have:
\begin{equation*}
G\theta_t|T(\ell)\equiv2G\theta_t\pmod{p^2},
\end{equation*}
where
\begin{equation*}
G\theta_t|T(\ell)=\sum_{n\in\mathcal{N}}c'(n)q^{n},
\end{equation*}
such that $c'(n)=c(n\ell)+\psi(\ell)\ell^{p-1}c(n/\ell)$ and
$c(n/\ell)=0$ whenever $\gcd{(n,\ell)}=1$. Hence, for every prime to
$\ell$ positive integer $n$, we have
\begin{equation*}
c(n\ell)\equiv2c(n)\pmod{p^2}.
\end{equation*}
Recall that
$\epsilon_{d}*\sigma_{2i-1}(m)=\sum_{r|m}\mu(r)\chi_{(-1)^{i}d}(r)r^{i-1}\sigma_{2i-1}(m/r)$.
We make the following hypothesis
for Serre's primes for the modular form $G\theta_t$.\\

$(H'_p)$:\; There exist a fundamental discriminant $t$ and a prime
number $\ell\in\mathcal{S}_t$ such that $\ell=ta^{2}+b^{2}$
and $\epsilon*\sigma_{2p-3}(b)\not\equiv1\pmod{p^{2}}$.\\

We will see that $(H_p)$ is equivalent to $(H'_p)$ if $b$ is a prime
number. Under the hypothesis $(H'_p)$ we obtain the following
result:

\begin{theo}\label{Theorem p}
Let $\mathcal{L}=\{\ell_1,...,\ell_s\}$ be a finite set of odd
primes. Let $p\geq5$ be a prime number. Assume that hypotheses
$(H'_p)$ is satisfied for some prime number $\ell$. There is a real
quadratic $p$-rational field $\mathbb{Q}(\sqrt{d})$ for some
fundamental discriminant $d<\ell$ such that $d\in\mathcal{N}$, which
means that $(\frac{d}{\ell_k})=1$ for every $\ell_k\in\mathcal{L}$.
\end{theo}

\noindent{Proof.} For every prime number $\ell$ satisfying Serre's
congruence for $G\theta_t$, we have
\begin{equation*}
c(\ell)\equiv2c(1)\pmod{p^2},
\end{equation*}
where $c(1)=p\zeta(2-p)$ is a $p$-adic unit. Write $c(\ell)$ as
follows
\begin{equation}\label{c{n}}
c(\ell)=A(\ell)+B(\ell),
\end{equation}
where
\begin{equation*}
A(\ell)=\sum_{\begin{array}{c}
             \ell=tx^2+y \\
             y\neq\square\\
             y\in\mathcal{N}
           \end{array}
}\alpha(x)h_1(p-1,y)
\end{equation*}
and
\begin{equation*}
B(\ell)=\sum_{\begin{array}{c}
             \ell=tx^2+y^2 \\
             y\in\mathcal{N}
           \end{array}}\alpha(x)h_1(p-1,y^2).
\end{equation*}
Then we have the following congruence
\begin{equation*}
A(\ell)\equiv2c(1)-B(\ell)\pmod{p^2}.
\end{equation*}
Moreover, $\ell\equiv1\pmod{4t}$ implies that the equation
\begin{equation*}
\ell=tx^2+y^2
\end{equation*}
has a unique solution $(a,b)$, see e.g., \cite[Chapter 1, page
31]{Cox}. Hence the sum $B(\ell)$ is non-trivial and
$$B(\ell)=2p\zeta(2-p)\epsilon*\sigma_{2p-3}(b).$$
We obtain the congruence
\begin{equation*}
A(\ell)\equiv2p\zeta(2-p)(1-\epsilon*\sigma_{2p-3}(b))\pmod{p^2}.
\end{equation*}
Since hypotheses $(H'_p)$ is satisfied for the prime number $p$,
there exist a square free integer $t$ and a prime number
$\ell\in\mathcal{S}_t$ such that $\ell=ta^{2}+b^{2}$, with $b$
satisfying the property
\begin{equation*}
\epsilon*\sigma_{2p-3}(b)\not\equiv1\pmod{p^2}.
\end{equation*}
This leads to the property
\begin{equation}\label{An}
A(\ell)\not\equiv0\pmod{p^2}.
\end{equation}
Moreover, we have by construction
\begin{equation*}
A(\ell)=\sum_{\begin{array}{c}
             \ell=tx^2+y \\
             y\neq\square
           \end{array}
}a(x)h_1(p-1,y)
\end{equation*}
such that each component of $A(\ell)$ is of $p$-adic valuation
$\geq1$. By (\ref{An}) there exist integers $x$ and $y$ such that
$\ell=tx^2+y$, $y$ is not a square and $v_p(\alpha(x)h_1(p-1,y))=1$.
Hence by definition of $h(p-1,y)$, there exists a fundamental
discriminant $d<\ell$ such that $y=dy_{1}^2$ and
$v_p(L(2-p,\chi_d))=0$. Then Proposition \ref{propo 2.1} and the
second statement of remark 3.2, gives that the field
$\mathbb{Q}(\sqrt{d})$ is $p$-rational and satisfies the
decomposition conditions.$\quad$ $\hfill\blacksquare$\vskip 6pt

The statement of Hypothesis $(H_p)$ in the introduction involves the
so called Wieferich primes defined as follows:

\begin{ef}
Let $a>1$ be an integer. A prime number $p\nmid a$ is said to be a
Wieferich prime of basis $a$ if
$$a^{p-1}-1\equiv0\pmod{p^2}.$$ Otherwise, the prime
$p$ is said to be non-Wieferich.
\end{ef}

Let $b=p_{1}^{n_1}...p_{s}^{n_s}$ be an integer. Remark that
$$\epsilon*\sigma_{2p-3}(b)=\prod_{i=1}^{s}\epsilon*\sigma_{2p-3}(p_{i}^{n_i}),$$
and
$$\epsilon*\sigma_{2p-3}(p_{i}^{n_i})=1+p_{i}^{p-2}\sigma_{2p-3}(p_{i}^{n_{i}-1})(p_{i}^{p-1}-1).$$
Suppose that the integer $b$ is a prime number, then
$$\epsilon*\sigma_{2p-3}(b)=1+b^{p-2}(b^{p-1}-1).$$ Hence
$\epsilon*\sigma_{2p-3}(b)\not\equiv1\pmod{p^{2}}$ precisely when
$p$ is a non-Wieferich prime of basis $b$. Hence $(H_p)$ is
equivalent to $(H'_p)$ in this case and Theorem \ref{wieferich} is
proved.

\begin{rem}
Let $i$ be an even integer and let $\mathcal{H}_i$ be the
Cohen-Eisenstein series of weight $i+\frac{1}{2}$. Under adequate
hypotheses $(H_{(p,i)})$, analogous arguments gives the existence of
real quadratic $(p,i)$-regular number fields with prescribed
arithmetic properties as in Theorem \ref{Theorem p}.
\end{rem}

\end{document}